\documentstyle[12pt,amssymbols]{article}  
\newcommand{\bq}{\begin{equation}}  
\newcommand{\eq}{\end{equation}}  
\newcommand{\bqa}{\begin{eqnarray}}  
\newcommand{\eqa}{\end{eqnarray}}  
\newcommand{\ra}{\rightarrow}

\def\ed{\end{document}}  
  
\def\ra{\rightarrow}  
  
\def\2pi{1\over 2\pi i}  

\def\~{\tilde}  
\def\newline{\hfil\break}

\def\ra{\rightarrow}

\def\sq2{\sqrt{2}}  
\def\sqk2{\sqrt{2(k+2}}  
\def\sqk{\sqrt{k}}

\def\be{\begin{equation}}  
\def\ee{\end{equation}}  
\def\br{\begin{array}}  
\def\er{\end{array}}  
\def\bea{\begin{eqnarray}}  
\def\eea{\end{eqnarray}}  
\def\ba{\begin{equation}\begin{array}}  
\def\ea{\end{array}\end{equation}}  
\def\bac{\begin{equation}\begin{array}{rll}}

\newcommand{\LS}{\widehat{sl(2)}}

\def\Z{{\Bbb Z}}

\begin{document}  
\rightline{ITP-SB-96-03}  
\rightline{December, 1995}  
\vbox{\vspace{-10mm}}  
\vspace{1.0truecm}  
\begin{center}  
{\LARGE \bf   
Bosonized Wakimoto construction in the principal gradation
 }\\[8mm]  
{\large A.H. Bougourzi }\\  
[6mm]{\it 
 ITP\\  
SUNY at Stony Brook\\  
Stony Brook, NY 11794}\\[20mm]  
  
\end{center}  
\vspace{1.0truecm}  
\begin{abstract}  
It is well known that the bosonized version of the Wakimoto 
construction allows the explicit realization of any affine
algebra $\widehat{g}$, with arbitrary level $k$ in the 
homogeneous gradation, 
in terms of $dim(g)$ free bosonic fields.
In this paper, we show in the case of 
the simplest affine algebra $\widehat{sl(2)}$,
that the bosonized 
Wakimoto realization can  be extended to the 
principal gradation only when $k$ is  equal 
to the critical level, i.e., -2.
In this case, this construction can be achieved in terms of
arbitrary number (larger than 1) of free bosonic fields. 
\end{abstract}  
\newpage  
\section{Introduction}  
  
It is well known that  affine algebras in their principal
gradation have the simplest  highest weight 
representations \cite{BeSe86}. 
In a recent paper \cite{Aual96},  the usefulness
of the principal gradation in the 8-vertex  model has 
been discussed. 
In the 
homogeneous gradation, it is well established that an explicit 
realization of affine algebras and their highest weight 
modules on which the central term acts with an arbitrary integer
$k$ 
is given in terms of free bosonic field (whose modes
satisfy Heisenberg algebras).
We refer to the latter realization 
as the bosonized form of the Wakimoto
construction \cite{Wak86}. The main feature of this construction
is that, for an arbitrary affine algebra $\widehat{g}$, it 
requires $dim(g)$ free bosonic fields for generic level $k$. 
Here $g$ is the corresponding finite Lie algebra with grade zero.
In particular, in the case of $\widehat{sl(2)}$, one needs
 three
free bosonic fields for generic $k$, except 
for $k=-2$ (the critical
level) and $k=4$ where only two fields are necessary
\cite{Bou92}. 
 
So far,  affine algebras in the principal gradation have been
realized only in terms of a combination of free bosonic fields
and parafermionic fields \cite{LeWi84}. 
However, the latter fields 
 have  such complicated generalized
commutation relations that they have not 
been useful so far in the practical computation of  
any physical quantity of models with affine symmetry.
The main purpose of this paper is then to find 
 a realization of $\LS$ (i.e., the simplest case) in the
principal gradation   
 in terms of just  free bosonic fields.
More specifically, 
we write the commutation relations as a current
algebra,   with defining relations given as operator product 
expansions (OPE's). Then, we briefly review the Lepowsky-Wilson
realization, which is valid only when $k=1$,
 in terms of a single bosonic field \cite{LeWi78}. 
As in the homogeneous 
case, we adapt this construction to generic level $k$ by
considering the most general expression for these currents in 
terms of a set of additional fields and a set of
 unknown parameters.

Consisenty conditions 
with the above OPE's lead to  a 
set
of constraints on these free parameters. Solving these
constraints, we show that, unlike in the homogeneous gradation, 
only the case with the 
critical level $k=-2$ survives this Wakimoto type of 
construction. Moreover, this can be acheived in terms of an
arbitrary number (larger that 1) of free bosonic fields.

\section{The current algebra $\LS$}

The algebra $\LS$ in the principal gradation is generated
by $\{k, x_m, a_n; m\in\Z; n\in 2\Z+1\}$, with defining relations
\cite{LeWi78,Lep85}:
\bac
{[k, a_n]}&=&{[k,x_n]}=0,\\
{[a_n,a_m]}&=&nk\delta_{n+m,0},\quad {[a_n,x_m]}=2x_{n+m},\\
{[x_n,x_m]}&=&2(-1)^{n+1}a_{n+m}+(-1)^n kn\delta_{n+m,0},
\label{eq1}
\ea
where $n\in\Z$ and it is understood that $a_{2n}=0$. From the
above relations, $k$ is then a central element
and acts as a scalar denoted by the same letter $k$ on the
$\LS$ highest weight modules. Let us introduce the 
following currents
(generating functions):
\bac
a(z)&=&\sum_{n\in 2\Z +1}a_nz^{-n-1},\\
x(z)&=&\sum_{n\in\Z}x_n z^{-n-1},
\label{eq2}
\ea
with $z$ being a complex variable.

The algebra (\ref{eq1}) can then be written as a current
algebra generated by the currents $a(z)$ and $x(z)$,
with the following defining operator product
expansions (OPE's):
\bea
a(z)a(w)&=&{k(z^2+w^2)\over (z^2-w^2)^2}+{\rm regular\,\, terms
\,\, as\,\,}
z\ra \pm w,\quad |z|>|w|,
\label{eq3}\\
a(z)x(w)&=&{2wx(w)\over z^2-w^2}+{\rm regular\,\, terms\,\, as
\,\,}
z\ra \pm w,\quad |z|>|w|,
\label{eq4}\\
x(z)x(w)&=&-{2a(w)\over z+w}-{k\over (z+w)^2}+
{\rm regular\,\, terms\,\, as\,\,}
z\ra  -w,\quad |z|>|w|.
\label{eq5}
\eea
As we will see in the sequel, 
the above current algebra makes the realization of $\LS$ 
in terms of free bosonic fields more apparent. 

\section{Wakimoto realization of $\LS$}

The purpose of this section is to solve, for arbitrary $k$,
 the OPE's (\ref{eq3})-(\ref{eq5}) for the currents $a(z)$ 
and $x(z)$ in terms of 
  free bosonic fields. 
For this, we  first recall the Lepowsky-Wilson realization
of $\LS$ with $k=1$ in terms of a 
single bosonic field $\phi_1(z)$. This field   turns out 
to be related to $a(z)$ itself as:
\be
a(z)=i\partial\phi_1(z),
\ee
that is,
\be
\phi_1(z)=i\sum_{n\in 2\Z+1}{a_n\over n} z^{-n}.
\label{eq51}\ee
Henceforth we refer to any generating function of this type 
such that its modes satisfy a Heisenberg algerba as
a free bosonic field.

\subsection {Lepowsky-Wilson realization of $\LS$}

As mentioned earlier,  this realization is 
valid just for $k=1$. In this case,
the OPE of $\phi_1(z)$ with itself takes the form:
\be
\phi_1(z)\phi_1(w)=<\phi_1(z)\phi_1(w)>+:\phi_1(z)\phi_1(w):,
\label{eq6}
\ee
with the symbol :: denoting the normal ordering, .i.e., 
the
creation modes $a_n$, $n<0$, are always placed to the left
of the annihilation modes $a_n$, $n>0$, and the vacuum-to-vacuum 
expectation value $<\phi_1(z)\phi_1(w)>$ 
being obtained from (\ref{eq1}) as:
\be
<\phi_1(z)\phi_1(w)>=-{1\over 2}\ln{{z-w\over z+w}}.
\label{eq7}
\ee
Relations (\ref{eq6}) and (\ref{eq7}) yield the following 
OPE's among vertex operators of types $i\alpha \partial 
\phi_1(z)$ and $\exp{i\beta\phi_1(z)}$, 
where $\partial$ denotes
partial derivative with respect to $z$, and $\alpha$ and $\beta$
are arbitrary complex parameters:
\bea
i\alpha\partial \phi_1(z) i\alpha^\prime \phi_1(w)&=&
-\alpha\alpha^\prime \partial_z\partial_w<\phi_1(z)\phi_1(w)>
-\alpha\alpha^\prime:\partial_z \phi_1(z) \partial_w 
\phi_1(w):,\nonumber\\
&=&{\alpha\alpha^\prime (z^2+w^2)\over (z^2-w^2)^2}
-\alpha\alpha^\prime:\partial_z \phi_1(z) \partial_w \phi_1(w):,
\label{eq8}\\
i\alpha\partial \phi_1(z) \exp{i\beta \phi_1(w)}&=&
-\alpha\beta\partial_z<\phi_1(z)\phi_1(w)>:
\exp{i\beta \phi_1(w)}:,\nonumber\\
&=&{\alpha\beta w\over z^2-w^2}:\exp{i\beta \phi_1(w)}:,
\label{eq9}\\
\exp{i\beta \phi_1(z)}\exp{i\beta^\prime \phi_1(w)}&=&
\exp\{-\beta\beta^\prime<\phi_1(z)\phi_1(w)>\}
:\exp\{i\beta \phi_1(z)+i\beta^\prime \phi_1(w)\}:,
\nonumber\\
&=&({z-w\over z+w})^{\beta\beta^\prime\over 2}
 :\exp\{i\beta \phi_1(z)+i\beta^\prime \phi_1(w)\}:.
\label{eq10}
\eea

Using these relations, it is easy to verify  
 that the currents
$a(z)$ and $x(z)$ as realized in terms of $\phi_1(z)$ 
\cite{LeWi78,Lep85}: 
\bea
a(z)&=&i\partial\phi_1(z),
\label{eq11}\\
x(z)&=&\epsilon {z^{-1}\over 2}\exp{2i \phi_1(z)},\quad
\epsilon=\pm 1,
\label{eq12}
\eea
satisfy the current algebra (\ref{eq3})-(\ref{eq5}) with $k=1$.
This algebra has then two highest weight modules 
(basic modules) isomorphic to
two envelopping algebras $U(\hat{h}^-)$, with $\epsilon =1$ and
$\epsilon=-1$, respectively. Here $\hat{h}^-$ is the subalgebra 
generated by $\{a_{n};n<0\}$. In this sense, this realization is
canonical since the whole $\LS$ algebra and its modules are
constructed from its Heisenberg subalgebra.

\subsection{Wakimoto realization of $\LS$}

The purpose of this section is to extend the previous realization
to arbitrary level $k$ {\rm \`a la} Wakimoto. This means that by
analogy with the usual Wakimoto realization in the homogeneous
case, our goal is to realize the currents $a(z)$ and 
$x(z)$ in terms of 
an $n$-dimensional bosonic vector field $\phi(z)=(\phi_1(z),
\phi_2(z),\dots,\phi_n(z))$ with
\be
<\phi_i(z)\phi_j(w)>=
-{\delta_{ij}\over 2}\ln{{z-w\over z+w}},\quad i,j=1,\dots,n.
\label{eq13}\ee
More specifically, by 
Wakimoto realization of $\LS$ we mean that both 
currents $a(z)$ and $x(z)$ are realized in terms of $\phi(z)$
as follows:
\bac
a(z)&=&i \alpha\partial\phi(z),
\label{eq14}\\
x(z)&=&i\beta\partial\phi(z)\exp{i\gamma \phi(z)},
\label{eq15}
\ea
where 
\bac
\alpha&=&(\alpha_1,\alpha_2,\dots,\alpha_n),\\
\beta&=&(\beta_1,\beta_2,\dots,\beta_n),\\
\gamma&=&(\gamma_1,\gamma_2,\dots,\gamma_n)   
\label{eq16}
\ea
are complex parameters to be determined consistently with 
(\ref{eq3})-
(\ref{eq5}). We note here that the scalar product
between these $n$-dimensional vectors is Euclidean. 
First, let us consider the realization of 
$a(z)$. It is clear from the previous subsection and
in particular (\ref{eq8}) 
that (\ref{eq3}) is satisfied if we simply rescale $a(z)$ as
\be
a(z)=i\sqrt{k}\partial\phi_1(z),
\label{eq17}
\ee  
that is, $\alpha=(\sqrt{k},0,\dots,0)$. Moreover, due to 
(\ref{eq9}) and (\ref{eq17}),  
relation (\ref{eq4}) is satisfied if the vectors
$\beta$ and $\gamma$ have the following forms:
\bac
\beta&=&(0,\beta_2,\dots,\beta_n),\\
\gamma&=&({2\over \sqrt{k}},\gamma_2,\dots,\gamma_n).  
\label{eq18}
\ea
Therefore, the remaining paramaters $\beta_2,\dots,\beta_n$
and 
$\gamma_2,\dots,\gamma_n$  are to be determined from the last
consistency relation (\ref{eq5}). For this purpose, we derive
the following OPE of
$x(z)$ as given by (\ref{eq15}) with itself:
\bac
x(z)x(w)&=&-\{{zw((\beta\gamma)^2-2 \beta^2)\over (z^2-w^2)^2}-
{\beta^2\over (z+w)^2}+{i(\beta\gamma)(z\beta\partial_z\phi(z)
-w\beta\partial_w\phi(w))\over z^2-w^2},\\
&&+(\beta\partial_z\phi(z))(\beta\partial_w\phi(w))\}
({z-w\over z+w})^{\gamma^2/2}
:\exp{i\gamma(\phi(z)+\phi(w))}:.
\label{eq19}
\ea  
Thus, for this relation to reduce to (\ref{eq5}), we 
impose the following constraints on the free parameters:
\bea
(\beta\gamma)^2-2 \beta^2&=&0,
\label{eq20},\\
\gamma^2&=&0,
\label{eq21}
\eea 
in which case it simplifies to
\be
x(z)x(w)={\beta^2\over (z+w)^2}
-{i(\beta\gamma)(\beta\partial_w\phi(w))
-i\beta^2\gamma\partial_w\phi(w)\over z+w}+
{\rm regular\,\, terms\,\, as\,\,} z\ra \pm w.
\label{eq22}
\ee
This latter relation coincides with (\ref{eq5}) provided
that the
following constraints are also imposed:
\bea
\beta^2&=&-k,
\label{eq23}\\
(\beta\gamma)\beta- \beta^2\gamma&=&(2\sqrt{k},0\dots,0).
\label{eq24}
\eea
To each solution of the equations (\ref{eq20}), (\ref{eq21}), 
(\ref{eq23}) and (\ref{eq24}) for the unknown 
parameters $\beta_2,\dots,\beta_n$ and 
$\gamma_2,\dots\gamma_n$, we obtain a genuine Wakimoto
realization of $\LS$ in the principal gradation. 
Let us now solve these equations. Relation (\ref{eq24}) implies
that:
\be
\gamma_i=-{(\beta\gamma)\over k}\beta_i,\quad i=2,\dots,n.
\label{eq25}
\ee
Using this relation and substituting $\beta^2$ by $-k$ in 
(\ref{eq21}), we arrive at
\be
\beta\gamma=2\epsilon,\quad \epsilon=\pm 1,
\label{eq26}
\ee
which, after  substituting $\beta^2$ by $-k$ 
in (\ref{eq20}),  leads to the unexpected and 
important result $k=-2$. We refer to the latter value of 
$k$ as the 
critical level for reasons 
related to the Sugawara construction in the homogeneous
gradation. Thus, relation (\ref{eq25}) reduces to
\be
\gamma_i=\epsilon\beta_i,\quad i=2,\dots,n.
\ee
To recapitulate, for any positive integer $n\geq 2$, we obtain
a Wakimoto realization of $\LS$ if  the 
following conditions
are satisfied:
\bac
k&=&-2,\\
\beta_1&=&0,\\
\gamma_1&=&{2\over \sqrt{k}},\\
\gamma_i&=&\epsilon\beta_i,\quad i=2,\dots,n;\\
\sum_{i=2}^n\beta_i^2&=&2.
\ea
From this, 
it is clear that for $k=-2$ we obtain a family of Wakimoto
 realizations of $\LS$ parametrized by the free parameters
$n$, $\beta_i$, ($i=2,\dots,n$), and $\epsilon$. 
In particular, the simplest realization in this family is 
obtained by setting 
\bac
\beta_2&=&\gamma_2=\sqrt{2},\\
\beta_i&=&\gamma_i=0,\quad i=3,\dots,n.
\ea

As mentioned earlier, it is an unexpected result that, 
unlike in the homogeneous case, the
Wakimoto realization of $\LS$ in the 
principal gradation exists only for the critical level 
$k=-2$. 
It is thus  interesting  to examine whether such
existence of the Wakimoto realization just in the critical level
extends to the general case and in particular to
$\widehat{sl(n)}$, since
its Wakimoto realization in the homogeneous case is relatively
simple, whereas for other algebras it is either unexisting or
extremely complicated for general level $k$.  

\section{Aknowlegements}

We are very grateful to  ITP for a postdoctoral fellowship.

\pagebreak

\end{document}